\documentclass[twoside, 10pt]{article}

\usepackage{mathrsfs,amsfonts,amsmath}
\RequirePackage{ifthen,calc}
\usepackage{theorem}

\def\R{{\mathbb R}}

\def\E{{\mathbb E}}
\def\P{{\mathbb P}}

\def\N{{\mathbb N}}

 \catcode`@=11
 \def\@evenhead{\hbox to\textwidth{\footnotesize\rm\thepage \hfill
  {\it }}} 

 \def\@oddhead{\hbox to \textwidth{\footnotesize{\it
 } \hfill\thepage}}

 \renewcommand{\section}{\makeatletter
 \renewcommand{\@seccntformat}[1]{{\csname the##1\endcsname.}\hspace{0.45em}}
 \makeatother \@startsection
{section}
{1}
{0pt}
{\baselineskip}
{0.5\baselineskip}
{\normalsize\bfseries\mathversion{bold}}}

\renewcommand{\subsection}{\makeatletter
 \renewcommand{\@seccntformat}[1]{{\csname the##1\endcsname.}\hspace{0.45em}}
 \makeatother \@startsection
{subsection}
{1}
{0pt}
{\baselineskip}
{0.5\baselineskip}
{\normalsize\bfseries\mathversion{bold}}}

\catcode`@=12

\newtheorem{thm}{\noindent Theorem}[section]
\newtheorem{lem}{\noindent Lemma}[section]
\newtheorem{cor}{\noindent Corollary}[section]

\theoremheaderfont{\normalfont\bfseries}
\theorembodyfont{\slshape} {\theorembodyfont{\rmfamily}} {\theorembodyfont{\rmfamily}
\newtheorem{defn}{\noindent Definition}[section]}
{\theorembodyfont{\rmfamily} } {\theorembodyfont{\rmfamily}
}
{\theorembodyfont{\rmfamily} } {\theorembodyfont{\rmfamily}
\newtheorem{rem}{\noindent Remark}[section]}
{\theorembodyfont{\rmfamily} } {\theorembodyfont{\rmfamily} } {\theorembodyfont{\rmfamily}
}

 \def\beqlb{\begin{eqnarray}}\def\eeqlb{\end{eqnarray}}
 \def\beqnn{\begin{eqnarray*}}\def\eeqnn{\end{eqnarray*}}
 \numberwithin{equation}{section}
\def\2R{\mathbb{R}_+\times\mathbb{R}}
\def\qed{\hfill$\square$\smallskip}

\def\3R{\mathbb{R}_+\times\mathbb{R}_-}

\begin{document}
\title{\bf  Convergence in law to operator fractional Brownian motions}
\author{\small Hongshuai Dai \thanks{E-mail:mathdsh@gmail. com} \\\small College of Mathematics and Information Sciences, Guangxi
University, \\\small Nanning 530004, China.
}

\maketitle

\begin{abstract}
\noindent In this paper, we provide  two approximations in law of
operator fractional Brownian  motions.  One is
constructed by Poisson processes, and the other generalizes  a result of  Taqqu (1975).

\medskip
\noindent{\bf 2000 Mathematics subject classification:}  60F17; 60G15

\medskip
\noindent{\bf Keywords:} Operator fractional Brownian motion,
Poisson processes, vector-valued Gaussian sequence, weak convergence.

\end{abstract}

\bigskip

\section{Introduction}
Self-similar processes, first studied rigorously by Lamperti \cite{LL}
under the name ``semi-stable", are stochastic processes that are invariant in distribution
under suitable scaling of time and space.
There has been an
extensive literature on self-similar processes. We refer to Vervaat
\cite{VW} for general properties, to Samorodnitsky and Taqqu
\cite{ST94} [Chaps.7 and 8] for studies on Gaussian and stable
self-similar processes and random fields.

The definition of self-similarity has been extended to allow for scaling by linear operators on $\R^d$, and the corresponding processes are called operator self-similar (o.s.s) processes in the literature.
 See Laha and Rohatgi
\cite{LR1982}, Hudson and Mason \cite{HM1982}, and Sato \cite{ST1982}.  Various examples of operator self-similar Gaussian and non-Gaussian processes have been constructed and studied by Maejima and Mason \cite{MM1994} and Mason and Xiao \cite{MX1999}. The theory of operator self-similarity runs somewhat parallel to that of operator stable measures and is also related to that of operator scaling random fields. See, for instance, Meerschaert and Scheffler \cite{Meerschaert2001}, Bierm$\acute{e}$ et al.\cite{HMH} and the references therein.

Let $End (\R^{d}) $ be the set of linear operators on $\R^d$
(endomorphisms) and let $Aut (\R^d)$ be the set of invertible linear
operators (automorphisms) in $End (\R^d)$. For convenience, we will
not distinguish an operator $D\in End (\R^d)$ from its associated
matrix relative to the standard basis of $\R^d$.

In this paper, we will use the following definition of o.s.s. processes, which corresponds to that  of Sato \cite{ST1982}, but is stronger than that of Hudson and Mason \cite{HM1982}

An $\R^d-$ valued stochastic  process $\tilde{Y}=\{\tilde{Y}(t)\}$ is said to be operator self-similar if
it is stochastically continuous, and there exists a $D\in End\{\R^d\}$ such that for every $c>0$
\beqlb\label{s1-1}
\tilde{Y}(ct) \stackrel{d}{=}c^D \tilde{Y}(t),\;t\in\R,
\eeqlb
where $\stackrel{d}{=}$ denotes equality of all finite-dimensional distributions, and $$
c^D=\exp \big((\log c)D)\big)=\sum_{k=0}^{\infty} \frac{1}{k!} (\log c)^k D^k.
$$
Any matrix for which (\ref{s1-1}) holds is called an {\it exponent} of the o.s.s process $\tilde{Y}$.

We say that a process $\tilde{Y}=\{\tilde{Y}(t)\}$ has stationary increments (s.i.) if for every $b\in\R$
\beqlb\label{Rs1-2} \tilde{Y}(\cdot+ b) - \tilde{Y}(b)\stackrel{d}{=}\tilde{Y}(\cdot) -\tilde{Y}(0).
\eeqlb

One of  examples of operator self-similar processes is the operator
fractional Brownian motion (OFBM). OFBMs are mean-zero, o.s.s., Gaussian processes with stationary increments.  They are of interest in several areas and for reasons similar to those in the univariate case. For example, they are obtained and used in the
context of multivariate time series and long range dependence (see,
for example,  Chung \cite{Chung2002}, Davidson and de Jong \cite{Davidson2000},   Davidson and Hashimzade \cite{Davidson2008}, Dolado and Marmol \cite{Dolado2004},  Robinson \cite{Robinson2008}, and Marinucci and Robinson \cite{Marinucci2000}). They are also studied in problems related to, for example, queuing systems and large deviations (see  Delgado \cite{Delgado2007}, and Konstantopoulos and Lin\cite{Lin1996}). In particular,  Mason and Xiao \cite{MX1999} studied the sample properties of a particular class of OFBMs. Didier and Pipiras \cite{Didier2009,Didier2011} studied the basic properties of OFBMs, such as the time reversibility, the behavior of the spectral density around zero, and so on.

On the other hand, weak convergence to FBM processes has been studied extensively
since the works of Davydov \cite{D1970} and Taqqu \cite{T1975}. In recent years
many new results on approximations of FBM processes have been
established. For example, Enriquez \cite{E2004} showed that a FBM can be
approximated in law by appropriately normalized correlated random
walks.  Delgado and Jolis \cite{DJ2000} proved that the law of a FBM can be weakly approximated by the law of some processes constructed from a standard Poisson process.

Despite a growing interest in OFBMs, there is little work that studies  weak limit theorems for OFBMs. In this paper we will study weak limit theorems for OFBMs. We will extend two approximation   results for the one-dimensional  FBM to the multivariate case of OFBM. The rest of this paper is organized as follows. In Section 2, we recall some preliminaries  and present the main results of this paper.  Based on  Poisson processes, we prove  weak limit theorem for OFBMs in Section 3.  Based on a stationary sequence,   we  study weak convergence to OFBMs in Section 4.

\section{Preliminaries}

We first recall some facts we need later. Throughout this paper, we will use $\|x\|_E$ to denote the usual Euclidean norm of $x\in \R^d$. For $A\in End (\R^d)$, let
 $\left\|A\right\|=\max_{\|x\|_E=1}\|Ax\|_E$ denote the operator norm of $A$. It is easy to see that for $A,B\in End({\R^d})$,
 \beqlb\label{s2-1}
 \left\|AB\right\|\leq \left\|A\right\|\cdot\left\|B\right\|,
 \eeqlb
and for every $A=(A_{ij})_{d\times d}\in End(\R^d)$,
\beqlb\label{s2-2}
\max_{1\leq i,j\leq d}|A_{ij}|\leq \left\|A\right\|\leq d^{\frac{3}{2}} \max_{1\leq i,j\leq d}|A_{ij}|.
\eeqlb
Let $\sigma (A)$ be the collection of all eigenvalues of $A$. We denote
\beqlb\label{s2-3}
\lambda_A=\min\{Re\lambda: \lambda\in\sigma (A)\}\; \hbox{and}\;\Lambda_A=\max\{Re\lambda: \lambda\in\sigma(A)\}.
\eeqlb

As standard for the multivariate context, we furthermore assume that the OFBM is {\it proper}.  A random variable in $\R^d$ is said to be {\it proper} if the support of its distribution is not contained in a proper hyperplane of $\R^d$.

 Let $x'$ denote the transpose vector of $x\in\R^d$, and $B^*$ be the adjoint operator of $B$. Let  $x_{+}=\max\{x,0\}$ and $x_{-}=\max\{-x,0\}$.

Didier and Pipiras \cite{Didier2009}   established the integral representations of OFBMs in the spectral domain. The following lemma comes from Didier and Pipiras \cite{Didier2009}[Theorem 3.1].
\begin{lem}\label{s5-lem1}
Let $D$ be a linear operator on $\R^d$ with $0<\Lambda_D,\lambda_D<1$. Let $X=\{X(t)\}$ be an OFBM with o.s.s. exponent $D$. Then  $X$ admits the integral representation
\beqlb\label{s5-1}
X(t)\stackrel{d}{=}\int_{\R}\frac{e^{itx}-1}{ix} \Big(x_+^{-(D-\frac{I}{2})}A+x_-^{-(D-\frac{I}{2})}\bar{A}\Big)W(dx)
\eeqlb
for some linear  operator $A$ on $\mathbb{C}^d$. Here, $\bar{A}$ denotes the complex conjugate and
\beqnn
W(x):=W_1(x)+iW_2(x)
\eeqnn
denotes a complex-valued multivariate Brownian motion such that $W_1(-x)=W_1(x)$ and $W_2(-x)=-W_2(x)$,  $W_1(x)$ and $W_2(x)$ are independent, and the induced random measure $W(x)$ satisfies
$$
\E\Big[W(dx)W^*(dx)\Big]=dx,
$$
where $W^*$ is the adjoint operator of $W$.
\end{lem}

\begin{rem} For fixed $t\in\R$, let
\beqnn
F(x)&&=\frac{e^{itx}-1}{ix} \Big(x_+^{-(D-\frac{I}{2})}A+x_-^{-(D-\frac{I}{2})}\bar{A}\Big)
\\&&=F_1(x)+iF_2(x).
\eeqnn
 Since $F(x)=\bar{F}(-x)$, $F_1(x)=F_1(-x)$ and  $F_2(x)=-F_2(-x)$.
\end{rem}

Inspired by Samorodnitsky  and Taqqu \cite{ST94} [Chap. 7],  up to a multiplicative constant, we can rewrite $\{X(t)\}$ as follows.
\beqlb\label{s5-2}
X(t)\stackrel{d}{=}\int_{0}^{\infty}G_1(x,t)W_1(dx)+\int_{0}^{\infty}G_2(x,t)W_2(dx),
\eeqlb
where
\beqnn
G_1(x,t)=\frac{\sin tx}{x} x^{-(D-\frac{I}{2})} A_1 + \frac{\cos tx -1}{x}x^{-(D-\frac{I}{2})} A_2,
\eeqnn
\beqnn
G_2(x,t)=\frac{\sin tx}{x} x^{-(D-\frac{I}{2})} A_2 + \frac{1-\cos tx }{x}x^{-(D-\frac{I}{2})} A_1,
\eeqnn
and$$
 A=A_1+ i A_2.
$$
\begin{rem}\label{s2-R-rem1}
Using the same method as Mason and Xiao \cite{MX1999}, we can get that for $i=1,2$,
\beqnn
\int_{\R_+}\big\|G_i(u,t)\big\|^2du <\infty.
\eeqnn
\end{rem}

In this paper, we study weak limit theorems for the OFBM $X$. We first recall some results about weak convergence to FBMs. Stroock \cite{S82} showed that
\begin{lem}\label{s2-lem3}
Let $\{N_{n}(t),\; n=1,2\cdots\}$ be a sequence of
 Poisson processes with intensity $n$, and construct the
continuous processes $\{\tilde{Y}_{n}(t)\}$  by
\beqlb\label{s2-6} \tilde{Y}_n(t)=\sqrt{n}\int_{0}^{t}(-1)^{N_{n}(s)}ds, \;
0\leq t\leq 1. \eeqlb Then the laws of $\{\tilde{Y}_n\}$ converge weakly,
 in  $\mathcal{C}([0,\;1])$, to the law
of a Brownian motion, as $n\to\infty$.
\end{lem}

Using the same method as Delgado and Jolis \cite{DJ2000}, one can easily get
\begin{cor}\label{s2-cor1}
Let $N_n(t)$ and $\hat{N}_n(t)$ be two independent Poisson processes with intensity $n$, and construct the processes $\{\hat{Y}_n(t)\}$  by
\beqlb\label{s2-cor2}
\hat{Y}_n(t)=\sqrt{n}\int_{\R_+} \frac{(1-\cos tx)}{x^{H+1/2}} (-1)^{N_n(x)}dx+\sqrt{n}\int_{\R_+} \frac{\sin tx}{x^{H+1/2}} (-1)^{\hat{N}_n(x)}dx,
\eeqlb
where $H\in(0,\;1)$.
 Then the laws of $\{\hat{Y}_n(t),n=1,2,\cdots\}_{t\in[0,\;1]}$ converge weakly in  $\mathcal{C}([0,\;1])$, to the law
of a fractional Brownian motion with index $H\in(0,\;1)$, as
$n\to\infty$.
\end{cor}

Inspired by Corollary \ref{s2-cor1}, we want to show that the OFBM $X$ given by (\ref{s5-2}) can also be approximated by a sequence of processes similar to (\ref{s2-cor2}).

Let $\theta_{n}^0(t)=\sqrt{n}(-1)^{N_{n}(t)}$ and
\beqlb\label{RS2-1}\theta_n(t)=\Big(\theta_{n}^1(t),\cdots,\theta_{n}^d(t)\Big)',\eeqlb where $\theta_{n}^i(t), i\in \{1,\cdots,d\}$ are independent copies of $\theta_{n}^0(t)$. Furthermore we define
\beqlb\label{s2-7}
\int_{0}^{t}\theta_n(s)ds=\Big(\int_{0}^{t}\theta_{n}^1(s)ds,\cdots,\int_{0}^{t}\theta_{n}^d(s)ds\Big)'.
\eeqlb

Inspired by Corollary \ref{s2-cor1} and the construction of $\{\hat{Y}_n\}$, we define the sequence $\{X_n(t)\}_{n\in\N}$ as
\beqlb\label{defX}X_n(t)=\int_{0}^\infty G_1(x,t)\theta_n(x)dx+ \int_{0}^\infty G_2(x,t)\hat{\theta}_n(x)dx,
 \eeqlb
where $\hat{\theta}_n(x)$ is an independent copy of $\theta_n(x)$. Then
 \begin{thm}\label{s2-thm1}
 The laws of $\{X_{n}(t),n=1,2,\cdots\}_{t\in[0,\;1]}$ given by (\ref{defX}) in
$\mathcal{C}^d[0,\;1]$ converge weakly to the law of the OFBM $X$ given
by (\ref{s5-2}), as $n\to\infty$, where $\mathcal{C}^d[0,\;1]=\mathcal{C}\big([0,\;1],\;\R^d\big)$.
 \end{thm}

On the other hand, Taqqu \cite{T1975} showed that a FBM can be approximated in law by normalized partial sums of stationary random variables.  Since  OFBMs are multivariate extensions of FBMs, it is interesting to extend Taqqu's \cite{T1975} result to the OFBM case.

Before  we state our result, we first introduce the following notation.
$$\{A(n)\}=\Big\{\big(A_{ij}(n)\big)_{d\times d}\Big\}\in End (\R^d)$$ and $$\{B(n)\}=\Big\{\big(B_{ij}(n)\big)_{d\times d}\Big\}\in End (\R^d)$$ are asymptotically equivalent,as $n\to\infty$, if  for any  $i,j\in\{1,\cdots,d\}$,  one of the following cases holds:
  \begin{itemize}
\item[(i)] There exists $N_0\in\N$ such that for all $n\geq N_0,$
$$B_{ij}(n)\neq 0 \;\textrm{and}\; \lim_{n\to\infty} A_{ij}(n) /B_{ij}(n) =1. $$
\item[(ii)]There exists $N_1\in\N$ such that for all $n\geq N_1,$
$$ B_{ij}(n)= 0\;\textrm{and} \;A_{ij}(n)=0.$$
\end{itemize}
We denote this as $A(n)\sim B(n)$, as $n\to\infty$.

Recall that a process $\tilde{Y}=\{\tilde{Y}(t)\}$ is time reversible if it  satisfies
\beqlb\label{s2-R-1}
\tilde{Y}(t)\stackrel{d}{=}\tilde{Y}(-t),\;t\in\R.
\eeqlb The author refers the readers to Didier and Pipiras \cite{Didier2009} for time reversible Gaussian processes with stationary increments.

Due to the following lemma, we only can extend Taqqu's result to the time reversible case.
\begin{lem}\label{Rs5-lem1}
Let $ \{\hat{Z}_i,i=1,2,\cdots\}$  be a stationary proper mean-zero Gaussian sequence of  $\R^d$-valued vectors, and $$\hat{S}_n(t)=\sum_{i=1}^{\left\lfloor nt \right\rfloor} \hat{Z}_i.$$
If all the finite-dimensional distributions of $\{\hat{S}_n(t)\}$
converge to the corresponding finite-dimensional distributions of a mean-zero proper Gaussian process $Z=\{Z(t)\}$. Then $Z$ is  time reversible.
\end{lem}
{\it Proof:} Define
\beqlb\label{Rs5-2}
\E\big[Z(t)Z'(s)\big]=R_Z(t,s).
\eeqlb
Since
$$
\hat{S}_n(t)\stackrel{F.D}{\Rightarrow}Z(t),\textrm{as}\quad n\to\infty,
$$
where $\stackrel{F.D}{\Rightarrow}$ denotes convergence of the finite-dimensional distributions,
we get that \beqlb\label{Rs5-3}\lim_{n\to\infty}\E\big[\hat{S}_n(t)\hat{S}'_n(s)\big]=R_Z(t,s).\eeqlb
Since $\{\hat{Z}_i\}$ is a stationary Gaussian sequence, we can get
\beqlb\label{Rs5-4}
\E\big[\hat{S}_n(t)\hat{S}'_n(s)\big]=\E\big[\hat{S}_n(s)\hat{S}'_n(t)\big]
\eeqlb
It follows from (\ref{Rs5-2}) and (\ref{Rs5-4}) that
\beqlb\label{Rs5-5}
R_Z(t,s)=R_Z(s,t).
\eeqlb
Since $\{\hat{Z}_i\}$ is a stationary Gaussian sequence, and $Z$ is a Gaussian process, we get that $Z$ has stationary increments.
Therefore, by Proposition 5.1 in Didier and Pipiras \cite{Didier2009} and (\ref{Rs5-5}), we get that the lemma holds. \qed

\begin{rem}\label{s2-rem}
It follows from Proposition 5.1 in Didier and Pipiras \cite{Didier2009}  that if the OFBM $X$ is  time reversible, then
$$
\E\big[X(t)X'(s)\big]=\frac{1}{2}\Big[|t|^D \Gamma |t|^{D^*}+ |s|^D \Gamma |s|^{D^*}-|t-s|^D \Gamma |t-s|^{D^*}\Big].
$$
\end{rem}
Now we state our result as follows.
\begin{thm}\label{s2-thm2}
Let $\{Z_i,i=1,2,\cdots\}$ be  a stationary proper mean-zero  Gaussian sequence of $\R^d$-valued vectors. We define  $$r(i,j)=\E[ Z_iZ'_j]=\Big(r_{kq}(|i-j|)\Big)_{d\times d},$$ Suppose that
\beqlb\label{s2-R-3}
\sum_{i=1}^N\sum_{j=1}^N r(i,j)\sim K BN^D \Gamma N^{D^*}B^*, \;\textrm{as}\; N\to\infty,
\eeqlb
where $\Gamma=\E[ X(1)X'(1)]$, $B\in Aut(\R^d)$ and $K>0$ is  a positive number.
Then
\beqlb\label{s2-11}
Q_N(t)=d_N \sum_{i=1}^{\left\lfloor Nt\right\rfloor} Z_i,
\eeqlb
with $d_N\sim CN^{-D}B^-$, converges weakly, as $N\to\infty$ in $\mathcal{D}^d[0,\;1]$, up to a multiplicative matrix from the left, to the time reversible OFBM $X$ given by (\ref{s5-2}) with $A_2A_1^*=A_1A_2^*$, where $C\in Aut(\R^d)$, and $\mathcal{D}^d[0,1]=\mathcal{D}\big([0,\;1],\R^d\big)$.
\end{thm}

We will prove Theorem \ref{s2-thm2} in Section 4 .

To end this section, we give a technical lemma which comes from
Maejima and Mason \cite{MM1994}.
\begin{lem}\label{s2-lem1}
Let $D\in End(\R^d) $. If $\lambda_D>0$
and $r>0$, then for any $\delta>0$, there exist positive constants $K_1$ and $K_2$ such that
\beqlb\label{s2-4}
\left\|r^D\right\| \leq \begin{cases}K_1 r^{\lambda_D-\delta},  &\textrm{for all}\;  r\leq 1,
\\K_2 r^{\Lambda_D+\delta}, &\textrm{ for all}\; r\geq 1.
\end{cases}
\eeqlb
\end{lem}

In the rest of this paper, most of the estimates  contain unspecified constants. An unspecified positive and finite constant will be denoted by $K$, which may not be the same in each occurrence. Sometimes we shall emphasize the dependence of these constants upon parameters.
\section{Weak limit theorem based on  Poisson processes}
In this section, we prove Theorem \ref{s2-thm1}. We first  state some technical results.
In order to prove the  main result, we need a tightness criterion for vector-valued stochastic processes.
\begin{lem}\label{s4-lem1}
Let $\big\{Z_n(t)\big\}_{n\in\N}$ be a sequence of stochastic processes in $\mathcal{D}^d[0,\;1]$ satisfying:
\begin{itemize}
\item[(i)] for every $n\in\N$, $Z_n(0)=0$ a.s.;
\item[(ii)] there exist constants $K>0$, $\beta>0$,  $\alpha>1$ and an integer $N_0\in\N$ such that
\beqlb\label{s4-a1}
\E\bigg[\Big\|Z_n(t)-Z_n(s)\Big\|_E^\beta\bigg] \leq K(t-s)^\alpha, n\geq N_0.
\eeqlb
\end{itemize}
Then $\{Z_n(t)\}$ is tight in $\mathcal{D}^d[0,\;1]$.
\end{lem}

The proof of Lemma \ref{s4-lem1} is classical (see Billingsley \cite{B1968} [Chap.3] and Ethier and Kurtz \cite{EK86} [Chap.3]). Here we omit the proof.

Next, we show that $\{X_n(t)\}$ is tight in $\mathcal{C}^d[0,\;1]$. We need the following lemmas.
\begin{lem}\label{s3-lem1}
For any even $m\in\N$, there exists a constant $K(m,d)>0$ such that for any measurable function $f:\R_+\to End(\R^{d})$ with $\int_{\R_+}||f(u)||^2 du<\infty$,
\beqlb\label{s3-1}
\E\Big[\Big\|\int_{\R_+} f(u)\theta_n(u)du\Big\|_E^m\Big] \leq K(m,d)\Big(\int_{\R_+}\left\|f(u)\right\|^2 du\Big)^{\frac{m}{2}},
\eeqlb
where $\theta_n(\cdot)$ is given by (\ref{RS2-1}).
\end{lem}

{\it Proof:}
Let $f(u)=\big(f_{ij}(u)\big)_{d\times d}=\big(f^1(u),\cdots,f^d(u)\big)'$, where $f^i(u)=\big(f_{i1}(u),\cdots,f_{id}(u)\big)$. Then
\beqlb\label{s3-2}
\int_{\R_+}f(u)\theta_n(u)du=\Big(\int_{\R_+}f^1(u)\theta_n(u)du,\cdots,\int_{\R_+}f^d(u)\theta_n(u)du\Big)'.
\eeqlb
It follows from (\ref{s3-2}) that
\beqlb\label{s3-3}
\E\Big[\Big\|\int_{\R_+}f(u)\theta_n(u)du\Big\|_E^m\Big]\leq K(m)\sum_{i=1}^{d}\E\Big[\Big |\int_{\R_+}f^i(u)\theta_n(u)du\Big|^m\Big].
\eeqlb

For each $i\in \{1,\cdots,d\}$, we have
\beqlb\label{s3-4}
\E\Big[\Big|\int_{\R_+}f^i(u)\theta_n(u)du\Big|^m\Big]\leq K(m) \sum_{j=1}^d \E\Big[\Big|\int_{\R_+}f_{ij}(u)\theta_{n}^j(u)du\Big|^m\Big].
\eeqlb

Since $\int_{\R_+}||f(u)||^2du <\infty$, it follows from (\ref{s2-2}) that
\beqlb\label{s3-5}
\int_{\R_+}|f_{ij}(u)|^2du \leq \int_{\R_+}\|f(u)\|^2du <\infty.
\eeqlb

Using the same method as in Delgado and Jolis \cite{DJ2000}, we have
\beqlb\label{s3-6}
\E\Big[\big|\int_{\R_+}f_{ij}(u)\theta_{n}^j(u)du\big| \Big]^m \leq K(m) \Big( \int_{\R_+}|f_{ij}(u)|^2du \Big)^{\frac{m}{2}}.
\eeqlb

Combining (\ref{s2-2}), (\ref{s3-3}), (\ref{s3-4}) and (\ref{s3-6}), we get that the lemma holds.\qed

\begin{lem}\label{s3-lem2}
 Choose $\delta>0$, such that $\lambda_D-\delta>0$ and $\Lambda_D+\delta<1$, then for any $0\leq s<t \leq 1$ and  even $m\in\N$, we have
 \beqlb\label{s3-7}
\E\bigg[\Big\|X_n(t)-X_n(s)\Big\|_E^{m}\bigg]\leq K(m,\delta,d)(t-s)^{mH},
 \eeqlb
 where $H= \lambda_D-\delta.$
\end{lem}

{\it Proof:} Using the inequality $(a+b)^m \leq 2^{m-1}(a^m+b^m)$ for $a, b>0$, we get that there exists a constant $K(m)>0$ such that
\beqlb\label{s3-8}
\E\bigg[\Big\|X_n(t)-X_n(s)\Big\|_E^m\bigg]\leq &&K (m)\E\bigg[\Big\|\int_{\R_+}\Big(G_1(x,t)-G_1(x,s)\Big)\theta_n(x)dx\Big\|_E^m\bigg]\nonumber
\\&&\quad+K(m)\E\bigg[\Big\|\int_{\R_+}\big(G_2(x,t)-G_2(x,s)\Big)\hat{\theta}_n(x)dx\Big\|_E^m\bigg]\nonumber
\\&&= K(m)I_1(m)+K(m)I_2(m),
\eeqlb
where
$$
I_1(m)=\E\bigg[\Big\|\int_{\R_+}\Big(G_1(x,t)-G_1(x,s)\Big)\theta_n(x)dx\Big\|_E^m\bigg],
$$
and
$$
I_2(m)=\E\bigg[\Big\|\int_{\R_+}\big(G_2(x,t)-G_2(x,s)\Big)\hat{\theta}_n(x)dx\Big\|_E^m\bigg].
$$
Since $\hat{\theta}_n(x)$ is an independent copy of $\theta_n(x)$, we only deal with $I_1(m)$. $I_2(m)$ can be done in the same method.

 By Lemma \ref{s3-1} and Remark \ref{s2-R-rem1}, we get that
 \beqlb\label{s3-9}
 I_1(m)\leq K(m,d) \Big(\int_{\R_+}\big\|G_1(x,t)-G_1(x,s)\big\|^2dx\Big)^{\frac{m}{2}}.
 \eeqlb

On the other hand,
\beqlb\label{s3-R-1}
\big\|G_1(x,t)-G_1(x,s)\big\|^2&&=\Big\|\big(\sin tx-\sin sx\big)(\frac{1}{x})^{D+\frac{I}{2}}\nonumber
\\&&\qquad+\big(\cos tx-\cos sx\big)(\frac{1}{x})^{D+\frac{I}{2}}\Big\|^2\nonumber
\\&&\leq K\Big\|\big(\sin tx-\sin sx\big)(\frac{1}{x})^{D+\frac{I}{2}}\Big\|^2\nonumber
\\&&\qquad+K\Big\|\big(\cos tx-\cos sx\big)(\frac{1}{x})^{D+\frac{I}{2}}\Big\|^2.
\eeqlb

By (\ref{s3-9}) and (\ref{s3-R-1}),
\beqlb\label{s3-R-3}
I_1(m)\leq &&K(m)\bigg[\int_{\R_+}\Big[\big\|(\sin tx-\sin sx)(\frac{1}{x})^{D+\frac{I}{2}}\big\|^2\Big]dx\bigg]^{\frac{m}{2}}\nonumber
\\&&+K(m)\bigg[\int_{\R_+}\Big[\big\|(\cos tx-\cos sx)(\frac{1}{x})^{D+\frac{I}{2}}\big\|^2\Big]dx\bigg]^{\frac{m}{2}}.
\eeqlb

For the sake of conciseness, we define
\beqnn
U_1(m)=\bigg[\int_{\R_+}\Big[\big\|(\sin tx-\sin sx)(\frac{1}{x})^{D+\frac{I}{2}}\big\|^2\Big]dx\bigg]^{\frac{m}{2}}
\eeqnn
and
\beqnn
U_2(m)=\bigg[\int_{\R_+}\Big[\big\|(\cos tx-\cos sx)(\frac{1}{x})^{D+\frac{I}{2}}\big\|^2\Big]dx\bigg]^{\frac{m}{2}}.
\eeqnn
 It follows from Lemma \ref{s2-lem1} that
 \beqlb\label{s3-10}
 U_2(m)&&\leq K(m,\delta,d)\Big(\int_{0}^1\big(\cos tx-\cos sx)^{2}\big(\frac{1}{x}\big)^{2(\Lambda_D+\delta+\frac{1}{2})}dx\Big)^{\frac{m}{2}}\nonumber
 \\&&\quad\qquad+K(m,\delta,d)
 \Big(\int_{1}^\infty\big(\cos tx-\cos sx)^{2}\big(\frac{1}{x}\big)^{2(\lambda_D-\delta+\frac{1}{2})}dx\Big)^{\frac{m}{2}}.
 \eeqlb

 The first term  on the r.h.s. of (\ref{s3-10}) can be bounded from above by
 \beqlb\label{s3-11}
&&K(m,\delta,d) \Big( \int_{0}^{1}  \sin^2 ((t-s)x/2)\big(\frac{1}{x}\big)^{2(\Lambda_D+\delta+\frac{1}{2})}dx\Big)^\frac{m}{2}\nonumber
\\&&\quad\qquad\leq
K(m,\delta,d)(t-s)^{m(\Lambda_D+\delta)}\Big(\int_{0}^{\infty}\frac{\sin^2 x}{x^{2(\Lambda_D+\delta+\frac{1}{2})}}dx\Big)^\frac{m}{2}\nonumber
\\&&\quad\qquad\leq
K(m,\delta,d)(t-s)^{m(\Lambda_D+\delta)}\Big(\int_{0}^{1}\frac{\sin^2 x}{x^{2(\Lambda_D+\delta+\frac{1}{2})}}dx+\int_{1}^\infty\frac{1}{x^{2(\Lambda_D+\delta+\frac{1}{2})}}dx \Big)^\frac{m}{2}.
\eeqlb

The second term on the r.h.s. of (\ref{s3-10}) can be bounded from above by
\beqlb\label{s3-12}
 &&K(m,\delta,d)\Big( \int_{1}^{\infty}\sin^2 ((t-s)x/2)\big(\frac{1}{x}\big)^{2(\lambda_D-\delta+\frac{1}{2})}dx\Big)^\frac{m}{2}\nonumber
\\&&\quad\qquad\leq K(m,\delta,d)
(t-s)^{m(\lambda_D-\delta)}\Big(\int_{0}^{\infty}\frac{\sin^2 x}{x^{2(\lambda_D-\delta+\frac{1}{2})}}dx\Big)^\frac{m}{2}\nonumber
\\&&\quad\qquad\leq K(m,\delta,d)
(t-s)^{m(\lambda_D-\delta)}\Big(\int_{0}^{1}\frac{\sin^2 x}{x^{2(\lambda_D-\delta+\frac{1}{2})}}dx+\int_{1}^{\infty}\frac{1}{x^{2(\lambda_D-\delta+\frac{1}{2})}}dx\Big)^\frac{m}{2}.
 \eeqlb

It follows from (\ref{s3-10}) to (\ref{s3-12}) that
\beqlb\label{s3-R-4}
U_2(m)\leq K(m,\delta,d)(t-s)^{mH},
\eeqlb
where $H= \lambda_D-\delta.$

Using the same method as the proof of (\ref{s3-R-4}), we  get that
\beqlb\label{s3-R-2}
U_1(m)\leq K(m,\delta,d)(t-s)^{mH}.
\eeqlb

It follows from (\ref{s3-R-3}), (\ref{s3-R-4}) and (\ref{s3-R-2}) that
 \beqlb\label{s3-13}
 I_1(m)\leq K(m,\delta,d) (t-s)^{mH}.
 \eeqlb

Using the same method as  the proof of (\ref{s3-13}), we get
 \beqlb\label{s3-14}
 I_2(m)\leq K(m,\delta,d)(t-s)^{mH}.
 \eeqlb
 The lemma follows from (\ref{s3-8}), (\ref{s3-13}) and (\ref{s3-14}).
 \qed
 \begin{rem}\label{s3-rem-1}
 It follows from Ethier and Kurtz \cite{EK86} [Chap.3 Proposition 10.3] that Lemma \ref{s3-lem2} also shows that $\P\big(X_n\in \mathcal{C}^d[0,\;1]\big)=1.$
 \end{rem}

In the following, we will prove that $\{X_n(t),\;t\in[0,\;1]\}$
given  by (\ref{defX}) converges weakly in $\mathcal{C}^d[0,\;1]$
to the OFBM $X$ given by
(\ref{s5-2}), as $n\to\infty$.  In order to obtain it, we need the following lemma.

\begin{lem}\label{s2-lem5} The laws of $\big\{\int_0^t\theta_n(u)du\big\}$ converge weakly in $\mathcal{C}^d[0,\;1]$ to the law of  an $\R^d$-valued
 Brownian motion, as $n\to\infty$.
\end{lem}

{\it Proof:}  For every even $m\in\N$ and any $0\leq s<t\leq 1$,
\beqlb\label{s2-8}
\E\bigg[\Big\|\int_{0}^{t}\theta_n(u)du-\int_{0}^{s}\theta_n(u)du\Big\|_E^m\bigg]\leq K(m) \sum_{i=1}^{d}\E\bigg[\Big|\int_{s}^{t}\theta_{n}^i(u)du\Big|^m\bigg].
\eeqlb
By Lemma 2 in Delgado and Jolis \cite{DJ2000},   we get
\beqlb\label{s2-9}
\E\bigg[\Big\|\int_{0}^{t}\theta_n(u)du-\int_{0}^{s}\theta_n(u)du\Big\|_E^m\bigg]\leq K(m)(t-s)^{\frac{m}{2}}.
\eeqlb
By Lemma \ref{s4-lem1} and Proposition 10.3 in Ethier and Kurtz \cite{EK86} [Chap.3],
(\ref{s2-9}) implies the tightness of $\{\int_{0}^{t}\theta_n(u) du\}$.

On the other hand,   the convergence of all finite-dimensional distributions follows from Lemma \ref{s2-lem3} and the well-known Cram$\acute{e}$r-Wold device (Jurek and Mason \cite{JM1993} [Chap.1] and Whitt \cite{W2002} [Chap.3]).

Finally tightness plus convergence of the finite-dimensional distributions implies the weak convergence. \qed

The proof of Theorem \ref{s2-thm1} is as follows.

 {\it Proof of Theorem \ref{s2-thm1}: }
 We first prove the tightness of the laws of $\{X_n(t),\; t\in[0,\;1]\}$. From Lemma \ref{s3-lem2}, for every even $m\in\N$, which satisfies $mH>1$, inequality (\ref{s3-7}) holds.  It follows from Lemma \ref{s4-lem1} that $\{X_n\}$ is tight.

Now we proceed with the identification of the limit law.  It is sufficient to prove that for any $q\in\N$,  $a_1,\cdots,a_q \in\R$ and $t_1\cdots,t_q\in [0,\;1]$,
 \beqlb\label{s4-1}
 \sum_{m=1}^q a_m X_n(t_m)\Rightarrow \sum_{m=1}^q a_m X(t_m).
 \eeqlb

We  define for any $m\in\{1,\cdots,q\}$,
\beqnn
&&\int_{\R_+}G_1(t_m,x) \theta_n(x)d x=\Big(\tilde{X}_{n}^1(t_m),\cdots,\tilde{X}_{n}^d(t_m)\Big)',
 \\&&\int_{\R_+}G_2(t_m, x)\hat{\theta}_n(x)d x=\Big(\hat{X}_{n}^1(t_m),\cdots,\hat{X}_{n}^d(t_m)\Big)',
\eeqnn
and
 \beqnn
&&\int_{\R_+}G_1(t_m,x) dW_1(x)=\Big(\tilde{X}^{1}(t_m),\cdots,\tilde{X}^{d}(t_m)\Big)',
 \\&&\int_{\R_+}G_2(t_m,x) d W_2(x)=\Big(\hat{X}^{1}(t_m),\cdots,\hat{X}^{d}(t_m)\Big)'.
 \eeqnn

By the Cram$\acute{e}$r-Wold device, in order to prove (\ref{s4-1}) we only need to show that  for any $(b^1,\cdots,b^d)'\in \R^d$
\beqlb\label{s4-3}
 \sum_{j=1}^{d}\sum_{m=1}^{q}a_m b^j \tilde{X}_{n}^j(t_m)\Rightarrow \sum_{j=1}^{d}\sum_{m=1}^{q}a_m b^j \tilde{X}^{j}(t_m),
 \eeqlb
 and
\beqlb\label{s4-4}
 \sum_{j=1}^{d}\sum_{m=1}^{q}a_m b^j \hat{X}_n^j(t_m)\Rightarrow \sum_{j=1}^{d}\sum_{m=1}^{q}a_m b^j \hat{X}^{j}(t_m),
 \eeqlb
since $\theta_n(u)$ and $\hat{\theta}_n(u)$ are independent.
 Now we only prove that (\ref{s4-3}) holds. (\ref{s4-4}) can be done in the same way.

In order to simplify the notation, let
$$
W_1(x)=\big(W_1^1(x),W_1^2(x),\cdots,W_1^{d}(x)\big)',
$$
and
 $$
 f(u,x)=G_1(u,x) =\big(f_{ji}(u,x)\big)_{d\times d}=\big(f^1(u,x),\cdots,f^d(u,x)\big)',
 $$
 where $f^j(u,x)=\big(f_{j1}(u,x),\cdots,f_{jd}(u,x)\big)$.

Therefore, we can rewrite (\ref{s4-3}) as follows.
\beqlb\label{s4-5}
\sum_{i=1}^d H_n^i\Rightarrow \sum_{i=1}^d H^i,
\eeqlb
where
\beqnn
H_n^i=\sum_{j=1}^d\sum_{m=1}^q b^ja_m
\int_{\R_+}f_{ji}(t_m,x)\theta_{n}^i(x)dx=\int_{\R_+}F_i(x)\theta_{n}^i(x)dx,
\eeqnn

\beqnn
H^i=\sum_{j=1}^d\sum_{m=1}^q b^ja_m
\int_{\R_+}f_{ji}(t_m,x) dW_1^i(x)=\int_{\R_+}F_i(x) d W_1^i(x),
\eeqnn
and
 $$
F_i(x)=\sum_{j=1}^d\sum_{m=1}^q b^ja_m f_{ji}(t_m,x).
$$

In order to prove (\ref{s4-5}), it is sufficient to prove that for any $\xi\in\R$,
\beqlb\label{s4-b1}
\E\Big[\exp\Big[i\xi\sum_{k=1}^d H_n^k\Big]\Big]\to\E\Big[\exp\Big[i\xi\sum_{k=1}^d H^k\Big]\Big]
\eeqlb
as $n\to\infty$.

Since $\theta_{n}^k,k=1,\cdots,d$  are mutually independent, and so are $W_1^k, k=1,\cdots,d$,
\beqnn
\\&&\Bigg|\E\bigg[\exp\Big[i\xi\sum_{k=1}^d H_n^k\Big]\bigg]-\E\bigg[\exp\Big[i\xi\sum_{k=1}^d H^k\Big]\bigg]\Bigg|
\\&&\quad\qquad\qquad=\bigg| \prod_{k=1}^d \E\Big[\exp\big[i\xi H^k_n\big]\Big]- \prod_{k=1}^d \E\Big[\exp\big[i\xi H^k\big]\Big]\bigg|
\\&&\quad\qquad\qquad\leq \bigg|\E \Big[\exp\big[i\xi H^1_n\big]\Big]
-\E\Big[\exp\big[i\xi H^1\big]\Big]\bigg|
\\&&\quad\qquad\qquad\qquad+\bigg| \prod_{k=2}^d \E\Big[\exp\big[i\xi H^k_n\big]\Big]-\prod_{k=2}^d \E \Big[\exp\big[i\xi H^k\big]\Big]\bigg|.
\eeqnn
By induction on $k=2,\cdots,d$, we  get that in order to prove (\ref{s4-b1}) it is sufficient to  prove that for $k\in\{1,\cdots,d\}$,
\beqlb\label{s4-6}
 H_n^k\Rightarrow H^k,
 \eeqlb
as $n\to\infty$.

 By (\ref{s2-2}), we  get that for every $k\in\{1,\cdots,d\}$,
 \beqlb\label{s4-7}
 \int_{\R_+}F^2_k(u)du \leq K\sum_{j=1}^d\sum_{m=1}^q \int_{\R_+}  f^{2}_{jk}(u,t_m)du<\infty.
 \eeqlb

 Therefore,  the proof of (\ref{s4-6}) follows the lines of the proof of Theorem 2.1 in  Dai and Li \cite{DL2009}.

 Combining (\ref{s4-5}) and (\ref{s4-6}), we get that (\ref{s4-3}) holds. Similarly, we get that (\ref{s4-4}) holds.

By Theorem 7.8  in Ethier and Kurtz \cite{EK86}[Chap.3] and Remark \ref{s3-rem-1}, we  get that Theorem \ref{s2-thm1} holds. The proof is done.\qed
\section{Weak convergence based on a stationary sequence}
In this section, we will prove  Theorem \ref{s2-thm2}. We first show that $\{Q_N(t):t\in [0,\;1]\}$ is tight in $\mathcal{D}^d[0,\;1]$.
Before we prove  tightness, we  give a technical lemma.
\begin{lem}\label{s6-lmea1}
Let $\{A(n)\},\{B(n)\}, \{C(n)\}$ and $\{D(n)\}\in End\{\R^d\}$. If  $$ A(n)\sim C(n),\textrm{as}\; n\to\infty,$$ and
$$B(n)\sim D(n),\textrm{as}\; n\to\infty,$$
then
$$
A(n)B(n)\sim C(n)D(n),
$$
as $n\to\infty$.
\end{lem}
One can easily get that the lemma holds. Here we omit the proof.
\begin{lem}\label{s6-lem1}
For any $0\leq s<t\leq 1$ and even $m\in\N$,  there exists $N_1\in\N$ such that
\beqlb\label{s6-1}
\E\bigg[\Big\|Q_N(t)-Q_N(s)\Big\|_E^m\bigg] \leq K (t-s)^{m(\lambda_D-\delta)},\quad N\geq N_1,
\eeqlb
where  $\delta>0$ with $\lambda_D-\delta>0$, and $Q_N(t)$ is given by (\ref{s2-11}).
\end{lem}
{\it Proof:} For any $s\leq t\in[0,\;1]$, we have
\beqlb\label{s6-2}
\E\bigg[\Big\|Q_N(t)-Q_N(s)\Big\|_E^m\bigg]=\E\bigg[\Big\|d_N\sum_{i=1}^{\left\lfloor  N(t-s)\right\rfloor}Z_i\Big\|_E^m\bigg]=\E\bigg[\Big\|Q_N(t-s)\Big\|_E^m\bigg],
\eeqlb
since $\{Z_i\}$ is a stationary sequence.

In order to simplify the notation, let
\beqnn
d_N\sum_{i=1}^{\left\lfloor N(t-s)\right\rfloor}Z_i=(Q_{\left\lfloor N(t-s)\right\rfloor}^1,\cdots,Q_{\left\lfloor N(t-s)\right\rfloor}^d)'.
\eeqnn
Since $(a+b)^m\leq 2^{m-1}(a^m+b^m)$ for $a,b>0$, we  get that
\beqlb\label{s6-3}
\E\bigg[\Big\|Q_N(t)-Q_N(s)\Big\|_E^m\bigg]\leq K(m) \sum_{i=1}^d \E\bigg[\Big|Q_{\left\lfloor N(t-s)\right\rfloor}^i\Big|^m\bigg].
\eeqlb
Since $\{Z_i\}$ is an  $\R^d$-valued Gaussian sequence, $Q_{\left\lfloor N(t-s)\right\rfloor}^i$ is  Gaussian.  So, we  get that $\E\bigg[\Big|Q_{\left\lfloor N(t-s)\right\rfloor}^i\Big|^m\bigg]$ is proportional to $\Bigg(\E\bigg[\Big|Q_{\left\lfloor N(t-s)\right\rfloor}^i\Big|^2\bigg]\Bigg)^\frac{m}{2}$.

It follows from (\ref{s6-3}) that
\beqlb\label{s6-4}
\E\bigg[\Big\|Q_N(t)-Q_N(s)\Big\|_E^m\bigg]&&\leq K  \Bigg(\E\bigg[\Big\|Q_N(t)-Q_N(s)\Big\|_E^2\bigg]\Bigg)^{\frac{m}{2}}.
\eeqlb

Note that for an $\R^d$-valued random variable $Q=(Q^1,\cdots,Q^d)'$,   $\E[\|Q\|_E^2]$ equals the sum of diagonal entries of the correlation matrix. So it follows from  (\ref{s2-2}) that \beqlb\label{s6-5} \E\bigg[\Big\|Q_N(t)-Q_N(s)\Big\|_E^m\bigg]&&\leq K
\bigg\|\E\Big[ \big[Q_N(t)-Q_N(s)\big]\big[Q_N(t)-Q_N(s)\big]'\Big]\bigg\|^\frac{m}{2}.
\eeqlb

By (\ref{s2-R-3}) and
Proposition 2.2.2 in  Meerschaert and Scheffler
\cite{Meerschaert2001},
 \beqlb\label{s6-6} \lim _{N\to\infty}\E
[Q_N(t)-Q_N(s)][Q_N(t)-Q_N(s)]'= K  C(t-s)^D \Gamma (t-s)^{D^*}C^*. \eeqlb

It follows from (\ref{s6-5}) and (\ref{s6-6}) that
there exists $N_1\in\N$ such that for all $N\geq N_1$,
\beqlb\label{s6-a1}
\E\bigg[\Big\|Q_N(t)-Q_N(s)\Big\|_E^m\bigg]&&\leq K \left\|(t-s)^D \Gamma (t-s)^{D^*}\right\|^{\frac{m}{2}},
\eeqlb
since $\left\|\cdot\right\|$ is continuous.

By (\ref{s6-a1}) and Lemma \ref{s2-lem1}, we get that for all $N\geq N_1$,
\beqlb\label{s6-a2}
\E\bigg[\Big\|Q_N(t)-Q_N(s)\Big\|_E^m  \bigg]&& \leq K\left\|(t-s)^D \right\|^{\frac{m}{2}} \times \left\|(t-s)^{D^*}\right\|^{\frac{m}{2}}\nonumber
\\&&\leq K (t-s)^{m(\lambda_D-\delta)}.
\eeqlb
  The proof is completed. \qed

Now we prove Theorem \ref{s2-thm2}.

{\it Proof of Theorem \ref{s2-thm2}:}  It follows from Lemmas \ref{s4-lem1} and \ref{s6-lem1} that  $\{Q_N(t)\}$ is tight.

Now  we prove the convergence of all finite-dimensional distributions. For any $p\in\N$ and  $t_1,\cdots, t_p\in[0,\;1]$, $Q_N(t_1),\cdots,Q_N(t_p)$ are jointly Gaussian, since $\{Z_i\}$ is a Gaussian sequence.

Define:
$$
R_N(t_i,t_j)=\E[Q_N(t_i)Q'_N(t_j)].
$$

In order to simplify the notation, we define
\beqnn
S_N(t)=\sum_{i=1}^{\left\lfloor Nt\right\rfloor}Z_i.
\eeqnn
Since $\{Z_i\}$ is a stationary Gaussian sequence, we have that if $t_i>t_j$,
\beqlb\label{s6-8}
\E[Q_N(t_i)Q'_N(t_j)]=&&\frac{1}{2}d_N \Big[\E[S_N(t_i)S'_N(t_i)]\nonumber
\\&&+\E[S_N(t_j)S'_N(t_j)]-\E[S_N(t_i-t_j)S'_N(t_i-t_j)]\Big]d^*_N.
\eeqlb
By (\ref{s2-R-3}),
\beqlb\label{s6-9}
\E[S_N(t_j)S'_N(t_j)]=\sum_{i=1}^{\left\lfloor Nt_j \right\rfloor}\sum_{k=1}^{\left\lfloor Nt_j \right\rfloor}r(i,k)\sim
 K BN^{D}t_j^D \Gamma t_j^{D^*} N^{D^*}B^*,
\eeqlb
as $N\to\infty$.

So, we  get
\beqlb\label{s6-10}
\lim_{N\to\infty} R_N(t_i,t_j)=R(t_i,t_j),
\eeqlb
where $R(t_i,t_j)=
\frac{K}{2}C\Big[|t_i|^D \Gamma |t_i|^{D^*}+ |t_j|^D \Gamma |t_j|^{D^*}-|t_i-t_j|^D \Gamma |t_i-t_j|^{D^*}\Big]C^*.
$

It follows from Remark \ref{s2-rem} that if $X$ is a time reversible OFBM, then
\beqlb\label{s5-R-1}
\E\Big[CX(t)(CX(s))'\Big]=\frac{1}{2}C\Big[|t|^D \Gamma |t|^{D^*}+ |s|^D \Gamma |s|^{D^*}-|t-s|^D \Gamma |t-s|^{D^*}\Big]C^*.
\eeqlb
Since $Q_N(t_1),\cdots,Q_N(t_p)$ are jointly Gaussian, it follows from (\ref{s6-10}) and (\ref{s5-R-1}) that
\beqnn
\big(Q_N(t_1),\cdots,Q_N(t_p)\big)\Rightarrow \sqrt{K}\big(CX(t_1),\cdots,CX(t_p)\big).
\eeqnn

By Theorem 7.8  in Ethier and Kurtz \cite{EK86}[Chap.3], we get that the convergence of finite-dimensional distributions and  tightness ensure  weak convergence of $\{Q_N(t)\}$. The proof is completed. \qed

At the end of this paper, we consider the OFBM given by Mason and Xiao \cite{MX1999}.
Mason and Xiao \cite{MX1999} introduced the OFBM $Y=\{Y(t)\}$ as follows.
\begin{defn}\label{s2-R-defn1}
Let $D$ be a linear operator on $\R^d$ with $0<\lambda_D,\;\Lambda_D<1$. For any $t\in\R$, define
\beqlb\label{defn-MX}
Y(t)=\int_{\R_+}\big(1-\cos(tx)\big)\big(\frac{1}{x}\big)^{D+I/2} d \tilde{W} (x) +\int_{\R_+}\sin tx \cdot \big(\frac{1}{x}\big)^{D+I/2} d \hat{W}(x),
\eeqlb
where $\tilde{W}$ is a vector-valued Gaussian measure, and  $\hat{W}$ is an independent copy of $\tilde{W} $.
\end{defn}

By (\ref{s5-1}) and (\ref{defn-MX}), we see that, up to a multiplicative constant, the OFBM $Y$  and the OFBM $X$ with $A=I$ have equal distributions.

Define
$$Y_n(t)=\int_{\R_+}\big(1-\cos(tx)\big)\big(\frac{1}{x}\big)^{D+I/2} \theta_n(x)d x +\int_{\R_+}\sin tx \cdot \big(\frac{1}{x}\big)^{D+I/2}\hat{\theta}_n(x) d x,$$ where $\theta_n(x)$ is given by (\ref{RS2-1})  and $\hat{\theta}_n(x)$ is an independent copy of $\theta_n(x)$.
We  immediately have

\begin{cor}
The laws of $\{Y_{n}(t),t\in[0,\;1]\}$  converge weakly to the law of $Y$ given
by (\ref{defn-MX}) in
$\mathcal{C}^d[0,\;1]$, as $n\to\infty$.
\end{cor}

By (\ref{defn-MX}), we get that
\beqlb\label{s2-R-2}
\E[Y(t)Y'(s)]=\frac{1}{2}\Big[|t|^D \tilde{\Gamma} |t|^{D^*}+ |s|^D \tilde{\Gamma} |s|^{D^*}-|t-s|^D \tilde{\Gamma} |t-s|^{D^*}\Big],
\eeqlb
where $\tilde{\Gamma}=\E\big[Y(1)Y'(1)\big]$. One can get that $\E\big[Y(t)Y'(s)\big]=\E\big[Y(s)Y'(t)\big]$. Therefore, by Proposition 5.1 in Didier and Pipiras \cite{Didier2009} , the OFBM $Y$ given by (\ref{defn-MX}) is  time reversible. Therefore, from Theorem \ref{s2-thm2}, we  get that

\begin{cor}\label{s4-R-thm2}
Let $\{\tilde{Z}_i,i=1,2,\cdots\} $ be  a stationary proper mean-zero Gaussian sequence of  $\R^d$-valued vectors. We define $$\tilde{r}(i,j)=\E[ \tilde{Z}_i\tilde{Z}'_j]=\Big(\tilde{r}_{kq}(|i-j|)\Big)_{d\times d}.$$  Suppose that
\beqlb\label{s2-10}
\sum_{i=1}^N\sum_{j=1}^N \tilde{r}(i,j)\sim K BN^D \tilde{\Gamma} N^{D^*}B^*, \;\textrm{as}\; N\to\infty,
\eeqlb
where $B\in Aut(\R^d)$ and $K>0$ is  a positive number.
Then
\beqlb\label{s4-R-11}
\tilde{Q}_N(t)=d_N \sum_{i=1}^{\left\lfloor Nt\right\rfloor} \tilde{Z}_i,
\eeqlb
with $d_N\sim C N^{-D}B^{-1}$, converges weakly as $N\to\infty$ in $\mathcal{D}^d[0,\;1]$,  up to a multiplicative matrix from the left, to the OFBM $Y$ defined by (\ref{defn-MX}).
\end{cor}

In the following, we give a special case where conditions in Theorem \ref{s2-thm2} and Corollary \ref{s4-R-thm2} are satisfied.
\begin{cor}\label{cor1} Suppose that  $D\in Aut(\R^d)$ is a real diagonal matrix with positive diagonal entries $\lambda_1,\cdots,\lambda_d$. Let $\{Z_i,i=1,2,\cdots\}$ be a stationary proper mean-zero Gaussian sequence of $\R^d$-valued vectors. We define  $r(|i-j|)=\E [Z_{j} Z'_i]$. Suppose that
\begin{itemize}
\item[(i)] Case $\frac{1}{2}<\lambda_D<1$:
\beqlb\label{s4-R-7}
r(j)\sim K \tilde{\Gamma} j^{2D-2I} \;\textrm{as}\; j\to\infty \;\textrm{with}\; K>0.
\eeqlb
\item[(ii)]
Case $0<\Lambda_D<\frac{1}{2}$:
\beqlb\label{s4-R-6}
r(j)\sim K \tilde{\Gamma} j^{2D-2I} \;\textrm{as}\; j\to\infty \;\textrm{with}\; K<0,
\eeqlb
and
\beqlb\label{s4-R-10}
r(0)+2\sum_{j=0}^{\infty}r(j)=\mathbb{O},
\eeqlb
where $\mathbb{O}$ denotes the zero matrix.
\end{itemize}
Then
\beqlb\label{s6-12}
Q_N(t)=d_N \sum_{i=1}^{\left\lfloor Nt\right\rfloor} Z_i,
\eeqlb
with $d_N\sim N^{-D}C^{-1}$, converges weakly as $N\to\infty$, up to a multiplicative constant, to the OFBM $Y$  given by (\ref{defn-MX}), where $C\in Aut(\R^d)$ is a diagonal matrix with diagonal entries $$\frac{1}{\sqrt{|2\lambda_1-1|2\lambda_1}},\cdots,\frac{1}{\sqrt{|2\lambda_d-1|2\lambda_d}}.$$
\end{cor}

{\it Proof:}    We will proceed to prove this corollary in two steps.

{\bf Step 1}.  We first deal with the case that $\frac{1}{2}<\lambda_D<1$.  Since $D$ is a diagonal matrix, we  get
that $j^{2D-2I}$ is still a diagonal matrix.  By some simple calculations, we get that the diagonal entries of  $j^{2D-2I}$
are
$$
j^{2\lambda_1-2},\cdots,j^{2\lambda_d-2}.
$$
For every $i\in\{1,\cdots,d\}$, we get from  Feller \cite{Feller1971}[p.281] and Taqqu \cite{T1975} that
\beqlb\label{cor-p-1}
\sum_{k=1}^j k^{2\lambda_i-2}\sim (2\lambda_i-1)^{-1}j^{2\lambda_i-1},\quad\textrm{as}\quad j\to\infty,
\eeqlb
since $\frac{1}{2}<\lambda_D\leq \Lambda_D<1$.

So, by (\ref{cor-p-1}),
\beqlb\label{cor-p-2}
\sum_{k=1}^j k^{2D-2I}\sim C_0 j^{2D-I},\quad\textrm{as}\quad j\to\infty,
\eeqlb
where $C_0\in Aut(\R^d)$ is a diagonal matrix with diagonal entries
$$
(2\lambda_1-1)^{-1},\cdots,(2\lambda_d-1)^{-1}.
$$
On the other hand, we have
\beqlb\label{s6-14}
E_N=\sum_{i=1}^N\sum_{j=1}^N r(i,j)=r(0)+\sum_{j=1}^{N-1}\big\{r(0)+2\sum_{i=1}^jr(i)\big\}.
\eeqlb

Since $D$ is a diagonal matrix, $\tilde{\Gamma}$ is  a diagonal matrix.
Hence, we  get
from (\ref{s4-R-7}), (\ref{cor-p-2}) and (\ref{s6-14}) that,  as $N\to\infty$,
\beqlb\label{cor-p-3}
E_N\sim  2KC_0 \tilde{\Gamma} \sum_{j=1}^N  j^{2D-I}\sim 2KC_0 C_1\tilde{\Gamma} N^{2D},
\eeqlb
where $C_1\in Aut(\R^d)$ is a diagonal matrix with diagonal entries
$$
(2\lambda_1)^{-1},\cdots,(2\lambda_d)^{-1}.
$$

Since $C_0$ and $C_1$ are all diagonal matrices, we choose a diagonal matrix $C\in Aut(\R^d)$ such that $C^2=C_0 C_1$. Since $C,N^D$ and $\tilde{\Gamma}$ are all diagonal matrices,
\beqlb\label{s4-R-9}
E_N\sim 2KCN^D\tilde{\Gamma}N^{D^*}C^*.
\eeqlb

Using the same method as the proof of Theorem \ref{s2-thm2}, we  get that  $\{Q_N(t)\}$ converges weakly, up to a multiplicative constant, to $Y$, as $N\to\infty$.

{\bf Step 2}.  We deal with the case that $0<\Lambda_D<\frac{1}{2}$.  Corresponding to (\ref{s4-R-6}), we assume that $K<0$ in the following.  For every $i\in\{1,\cdots,d\}$, we  get from  Feller \cite{Feller1971}[p.281] and Taqqu \cite{T1975} that
\beqlb\label{cor-p-4}
\sum_{k=j}^\infty k^{2\lambda_i-2}\sim -(2\lambda_i-1)^{-1}j^{2\lambda_i-1},\quad\textrm{as}\quad j\to\infty,
\eeqlb
since $0<\Lambda_D<\frac{1}{2}$.
Then by (\ref{cor-p-4}),
\beqlb\label{cor-p-5}
\sum_{k=j}^\infty k^{2D-2I}\sim \hat{C}_0 j^{2D-I},\quad\textrm{as}\quad j\to\infty,
\eeqlb
where $\hat{C}_0\in Aut(\R^d)$ is a diagonal matrix with diagonal entries
$$
-(2\lambda_1-1)^{-1},\cdots,-(2\lambda_d-1)^{-1}.
$$

On the other hand, it follows from (\ref{s4-R-10}) that
\beqlb\label{cor-p-6}
r(0)+2\sum_{i=1}^jr(i)=-2\sum_{i=j+1}^\infty r(i).
\eeqlb
By (\ref{s4-R-6}), (\ref{cor-p-5}) and (\ref{cor-p-6}),
\beqlb\label{cor-p-7}
r(0)+2\sum_{i=1}^jr(i)\sim -2K\hat{C}_0\tilde{\Gamma} j^{2D-I}.
\eeqlb

By (\ref{cor-p-1}),  (\ref{s6-14}), (\ref{cor-p-4}) and (\ref{cor-p-7}), we get that as $N\to\infty$,
\beqlb\label{cor-p-8}
E_N\sim -2K\hat{C}_0\tilde{\Gamma}\sum_{j=1}^N j^{2D-I} \sim-2K\hat{C}_0C_1 \tilde{\Gamma} N^{2D}.
\eeqlb

Since both $\hat{C}_0$ and $C_1$  are  diagonal matrices, we  choose $\hat{C}\in Aut(\R^d)$ such that $(\hat{C})^2=\hat{C}_0 C_1$. Therefore
\beqlb\label{s4-R-8}
E_N\sim-2K\hat{C}N^D\tilde{\Gamma}N^{D^*}\hat{C}^*,
\eeqlb
as $N\to\infty$.
Using the same method as the proof of Theorem \ref{s2-thm2}, we  get that  $\{Q_N(t)\}$ converges weakly, up to a multiplicative constant, to $Y$, as $N\to\infty$.  The proof has been completed.
\qed

\medskip
\noindent{\bf Acknowledgements} The author thanks Professor Yimin Xiao, Michigan State
University, U.S.A., and Professor Yuqiang Li, East China Normal University, China,  for stimulating discussions. I also would like to thank the reviewer for helpful comments to improve this work. This work was supported by the Scientific Research Foundation of Guangxi University (NO: XBZ110398).

\end{document}